\documentclass[final,3p,times]{elsarticle}
\usepackage{mathrsfs}
\usepackage{amsfonts}
\usepackage{amsmath,amssymb}
\usepackage[all]{xy}
\usepackage{latexsym}
\usepackage{amsthm,a4wide,color}
\usepackage{amsmath,amscd,verbatim,enumerate}
\usepackage{hyperref}
\usepackage{graphicx}
\usepackage{extarrows}
\usepackage{lineno}
\biboptions{sort&compress}

\input diagxy

\theoremstyle{plain}
  \newtheorem{thm}{Theorem}[section]
   \newtheorem{theorem}[thm]{Theorem}
  \newtheorem{lemma}[thm]{Lemma}
  \newtheorem{proposition}[thm]{Proposition}
  
\theoremstyle{definition}
  \newtheorem{definition}[thm]{Definition}
  
  \newtheorem{example}[thm]{Example}
  \newtheorem{remark}[thm]{Remark}

\begin{document}

\newcommand{\oto}{{\to\hspace*{-3.1ex}{\circ}\hspace*{1.9ex}}}
\newcommand{\lam}{\lambda}
\newcommand{\da}{\downarrow}
\newcommand{\Da}{\Downarrow\!}
\newcommand{\D}{\Delta}
\newcommand{\ua}{\uparrow}
\newcommand{\ra}{\rightarrow}
\newcommand{\la}{\leftarrow}
\newcommand{\lra}{\longrightarrow}
\newcommand{\lla}{\longleftarrow}
\newcommand{\rat}{\!\rightarrowtail\!}
\newcommand{\up}{\upsilon}
\newcommand{\Up}{\Upsilon}
\newcommand{\ep}{\epsilon}
\newcommand{\ga}{\gamma}
\newcommand{\Ga}{\Gamma}
\newcommand{\Lam}{\Lambda}

\newcommand{\TF}{\mathbb{F}}
\newcommand{\TU}{\mathbb{U}}
\newcommand{\TW}{\mathbb{W}}
\newcommand{\TG}{\mathbb{G}}
\newcommand{\TB}{\mathbb{B}}
\newcommand{\TV}{\mathbb{V}}
\newcommand{\TT}{\mathcal{T}}
\newcommand{\TS}{\mathcal{S}}

\newcommand{\TC}{\mathbb{C}}
\newcommand{\CU}{\mathscr{U}}
\newcommand{\CW}{\mathscr{W}}
\newcommand{\CCF}{\mathscr{F}}
\newcommand{\CCG}{\mathscr{G}}
\newcommand{\CVV}{\mathscr{V}}
\newcommand{\CN}{{\mathcal{N}}}

\newcommand{\SU}{{\mathcal{U}}}
\newcommand{\SV}{{\mathcal{V}}}
\newcommand{\SF}{{\mathcal{F}}}
\newcommand{\SG}{{\mathcal{G}}}
\newcommand{\SB}{{\mathcal{B}}}
\newcommand{\ST}{{\mathcal{T}}}

\newcommand{\DV}{{\mathbf{V}}}

\newcommand{\DB}{{\mathbf{B}}}

\numberwithin{equation}{section}
\renewcommand{\theequation}{\thesection.\arabic{equation}}

\begin{frontmatter}
\title{The further study on the category of $\top$-convergence groups\tnoteref{F}}

\tnotetext[F]{This work is supported by National Natural Science Foundation of China (No. 12171220).}

\author{Lingqiang Li\corref{cor}}
\ead{lilingqiang0614@126.com}\cortext[cor]{Corresponding author}
\author{Qiu Jin}
\ead{jinqiu79@126.com}

\address{Department of Mathematics, Liaocheng University, Liaocheng 252059, P.R.China}

\begin{abstract} $\top$-convergence groups is a natural extension of lattice-valued topological groups, which is a newly introduced mathematical structure. In this paper, we will further explore the theory of $\top$-convergence groups. The main results include: (1) It possesses a novel characterization through the $\odot$-product of $\top$-filters, and it is localizable, meaning that each $\top$-convergence group is uniquely determined by the convergence at the identity element of the underlying group. (2) The definition of its subcategory, the topological $\top$-convergence groups, can be simplified by removing the topological condition (TT). (3) It exhibits uniformization, which means that each $\top$-convergence group can be reconstructed from a $\top$-uniformly convergent space. (4) It possesses a power object, indicating that it has good category properties.

\end{abstract}


\end{frontmatter}

\section{Introduction}

Topological group is an important mathematical structure, which is defined as a topological space over a group s.t. the group operations are continuous \cite{AA08}. As natural extensions, fuzzy topological groups are proposed by replacing topological spaces with fuzzy topological spaces. Due to the diversity of fuzzy topologies, there are many kinds of fuzzy topological groups.
\begin{itemize}

\item[$\diamond$] Foster \cite{DF79} defined fuzzy topological group by Lowen fuzzy topology \cite{RL76}, Ahsanullah \cite{TA88} introduced fuzzy neighborhood group through Lowen fuzzy neighborhood space \cite{RL82}, and Yan and Guo \cite{CY10} investigated $I$-fuzzy topological group via  $I$-fuzzy quasi-coincident neighborhood system \cite{JF04}.

\item[$\diamond$] Bayoumi \cite{FB05} defined $L$-topological group by Chang-Goguen $L$-topology, Bayoumi \cite{FB08} introduced Global $L$-neighborhood group via global $L$-neighborhood system, Zhang and Yan \cite{Yan12} investigated $L$-fuzzifying topological group through $L$-fuzzifying neighborhood system, and Zhao et al. \cite{DZ14} gave $(L,M)$-fuzzy topological group via  $(L,M)$-neighborhood system \cite{FS09}, where $L$ and $M$ are completely distributive complete lattices. Zhang and Fang \cite{HZ07}  discussed $I(L)$-topological groups and their level $L$-topological groups.

\item[$\diamond$] Mufarrij and Ahsanullah \cite{JA08} defined stratified $L$-topological group by  stratified $L$-neighborhood space \cite{U.H99}, Ahsanullah and J\"{a}ger introduced \cite{TA19} quantale-valued approach group through approach space, Shen \cite{SH92,SH94} proposed fuzzifying topological group via
$L$-fuzzyfing topological space, and Li and Jin \cite{LL21} investigated $\top$-neighborhood group
by $\top$-neighborhood space, where $L$ is a complete
residuated lattice (i.e., a commutative integral quantale)\footnote{Complete
residuated lattices include GL-monoid, BL-algebra, MV-algebra,
$\prod$-algebra, Heyting algebra, and left continuous t-norm, etc.
as special case and can be regarded as truth table of generalized
multi-valued logic.}. The category of $\top$-neighborhood groups is proven to be equivalent to that of strong $L$-topological groups.

\end{itemize}

Fuzzy convergence spaces are generalizations of fuzzy topological
spaces by relaxing the convergence of fuzzy topological spaces. The stratified $L$-generalized convergence spaces defined by $L$-filters \cite{G.J01} and $\top$-convergence spaces defined by $\top$-filters \cite{fang2017} are two well known  fuzzy convergence spaces. Note that topological stratified $L$-generalized convergence spaces  and  topological $\top$-convergence spaces respectively characterize two significant types lattice-valued topological spaces: stratified $L$-topological spaces \cite{G.J07} and strong $L$-topological spaces \cite{Yu17}. Nowadays, these two kinds of lattice-valued convergence spaces have been extensively studied and expanded, see  Fang et al. \cite{J.F10,JF21,JF23,Yu17}, J\"{a}ger \cite{G.J15,GJ19,G.J24,YY22}, Lai and Zhang \cite{HL17}, Li et al. \cite{LL12,LL18}, Pang et al. \cite{P17,P172,P18,LZ22}, Reid and Richardson \cite{GR,GR18}, Yao \cite{WY12} and Yue et al. \cite{Yue20,Yue21}.


Recently, lattice-valued convergence spaces have been utilized to supplant lattice-valued topological spaces for constructing a broader lattice-valued convergence group theory than the lattice-valued topology group theory.
\begin{itemize}

\item[$\diamond$] For a frame $L$, Ahsanullah and Mufarrij \cite{TA08} pioneered the concept of a lattice-valued convergence group, specifically a frame-valued stratified generalized convergence group, which is defined as a stratified $L$-generalized convergence space \cite{G.J01} over a group, with the condition that the group operations are continuous. They also provided an intriguing equivalent characterization of their lattice-valued convergence groups through the $\odot$-product of stratified $L$-filters. This characterization has proven useful in the analysis of properties within frame-valued stratified generalized convergence groups. Moreover, it was demonstrated that frame-valued stratified generalized convergence groups represent an extension of stratified $L$-topological groups, as defined by stratified $L$-neighborhood systems \cite{JA08}. Subsequently, Ahsanullah et al. \cite{TG142} explored enriched $L$-generalized convergence groups, where $L$ is a GL-monoid. Ahsanullah and J\"{a}ger \cite{TA17} introduced the stratified $LMN$-convergence tower group, involving frames $L$ and $M$ and a complete residuated lattice $N$. Ahsanullah and J\"{a}ger \cite{TA19} presented quantale-valued generalizations of approach groups. It's important to note that the lattice-valued convergence groups in \cite{TG142, TA17, TA19} are defined via the $\odot$-product of lattice-valued filters, rather than through the continuity of group operations. 

 \item[$\diamond$] For a complete residuated lattice $L$ , Zhang and Pang \cite{LZ23} introduced a novel type of lattice-valued convergence groups via $\top$-filters, denoted as $\top$-convergence groups. This is defined as a $\top$-convergence space \cite{fang2017} over a group with the property that group operations are continuous. They demonstrated that the $\top$-convergence groups extend the $\top$-neighborhood groups, and also generalize strong $L$-topological groups, considering the isomorphism between these two categories as shown in \cite{LL21}. Nonetheless, Zhang and Pang did not provide an equivalent characterization on their $\top$-convergence groups through the $\odot$-product of $\top$-filters, which may pose difficulties for subsequent research into $\top$-convergence groups.
\end{itemize}

After the comprehensive review, we have identified that two types of lattice-valued convergence groups have garnered more interest within the context of complete residual lattices. The first type is the enriched $L$-generalized convergence groups \cite{TG142}, characterized by the stratified $L$-filters.  The enriched $L$-generalized convergence groups are  the extensions of stratified $L$-topological groups proposed in \cite{JA08}. The second type is the $\top$-convergence groups \cite{LZ23}, defined by the $\top$-filters. The $\top$-convergence groups are the generalizations of strong $L$-topological groups given in \cite{LL21}. In the theory of (fuzzy) topological groups, homogeneity, localization, uniformization and category property constitute important research subjects. Following (fuzzy) topological group theory, the localization property of lattice-valued convergence groups should mean that lattice-valued convergence is uniquely determined by the convergence at the identity element of the underlying group. Moreover, the uniformization property should indicate that each lattice-valued convergence group can be reconstructed from a lattice-valued uniformly convergence space.
\begin{itemize}

\item[$\diamond$] For enriched $L$-generalized convergence groups: Ahsanullah et al. \cite{TG142} showed that the enriched $L$-generalized convergence groups possess homogeneity and constitute topological categories over groups. They further demonstrated that each enriched $L$-generalized convergence group induces an $L$-uniformly convergence space; however, they failed to prove that $L$-uniformly convergence space can, in turn, reproduce the enriched $L$-generalized convergence group. Consequently, the uniformization of enriched $L$-generalized convergence group remains incomplete. Additionally, the localization of enriched $L$-generalized convergence group has not been studied.  Ahsanullah and Mufarrij \cite{TA08} verified that framed-valued stratified generalized convergence groups possesses a power object; yet, their result cannot be generalized to the context of complete residuated lattice case.

\item[$\diamond$] For $\top$-convergence groups: Zhang and Pang \cite{LZ23} demonstrated that $\top$-convergence groups exhibit homogeneity and form a topological category over groups. The subcategory of topological $\top$-convergence groups within the category of $\top$-convergence groups is isomorphic to the category of topological $\top$-neighborhood groups. Furthermore, the subcategory of $\top$-limit groups within the category of $\top$-convergence groups possesses uniformization; that is, each $\top$-limit group can be reconstructed from a $\top$-uniformly limit space. However, the uniformization,  localization and power objects of $\top$-convergence groups have not yet been addressed.
\end{itemize} 

 The research status of two types of lattice-valued convergence groups is illustrated in Figure \ref{fig11} below.
\begin{figure}[!ht]
\centering
\includegraphics[width=14cm,height=5.5cm]{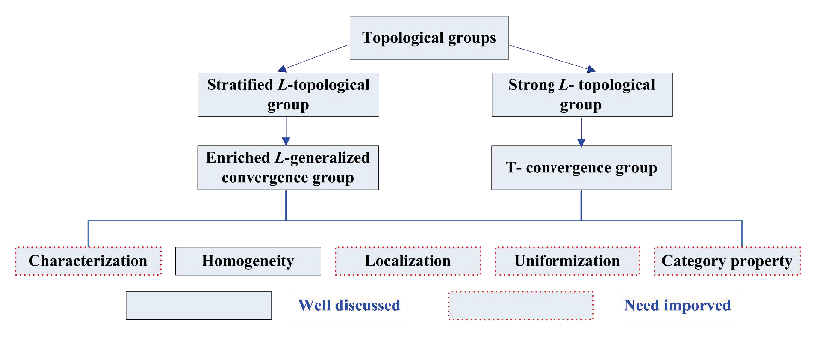}
\caption{ The research status on enriched $L$-generalized convergence groups and $\top$-convergence groups} \label{fig11}
\end{figure}

As illustrated in  Figure \ref{fig11}, there is a substantial requirement for improvement in the characterization, localization, uniformization and category property of complete residuated lattice-valued $\top$-convergence groups. In this paper, we will undertake additional investigations into these aforementioned issues.

The contents are arranged as follows. In Section 2, we recall some basic
concepts as preliminary. In Section 3, we characterize the $\top$-convergence group via the $\odot$-product of  $\top$-filters, and subsequently derive its localization. In Section 4, we introduce the concept of pretopological $\top$-convergence groups by eliminating the topological condition (TT) from topological $\top$-convergence groups as presented in \cite{LZ23}. We then demonstrate that the category of pretopological $\top$-convergence groups is equivalent to that of topological $\top$-convergence groups. Consequently, the topological condition (TT) in \cite{LZ23} is deemed redundant. In Section 5, we show that category of $\top$-convergence groups has uniformization, that is,  each $\top$-convergence group can be
reproduced by a $\top$-uniformly convergence space. In Section 6, we  prove that the category of $\top$-convergence groups has power object. In Section 6, we make a conclusion.

\begin{figure}[!ht]
\centering
\includegraphics[width=14cm,height=3cm]{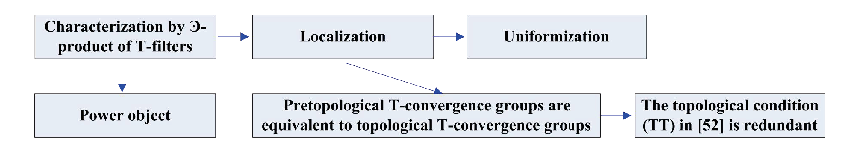}
\caption{The mind map of this paper} \label{fig44}
\end{figure}

The mind map of this paper can be summarized as shown in Figure \ref{fig44}.

\section{Preliminaries}

In this section, we recall some notions and notations for later use.

A {\it complete residuated lattice} \cite{H98} is a pair $(L,\ast)$ satisfying:

(1) $(L,\leq)$ is a complete lattice with  top
(resp., bottom) element $\top$ (resp., $\bot$),

(2) $\ast$ is a commutative and  associative binary operation on $L$ s.t. $\forall a\in L, \forall \{ b _j\}_{j\in J}\subseteq L$:
$$a \ast\bigvee\limits_{j\in J} b _j=\bigvee\limits_{j\in J}( a \ast  b _j),\ \ \ a\ast \top=a.$$

In this paper, we always assume that $L=(L,\ast)$ represent a complete residuated lattice.

The operation $\ast$ has a residuated operation $\rightarrow: L\times L\longrightarrow L$ determined by  $a\rightarrow b=\bigvee\{c\in L|a\ast c\leq b\}$. Their properties are summarized as below:

\begin{enumerate}
\item[{\rm (I1)}] $a \rightarrow b=\top\Longleftrightarrow a \leq b$,
\item[{\rm (I2)}] $a \ast(a \rightarrow b)\leq b$,
\item[{\rm (I3)}] $a \rightarrow(b\rightarrow c )=(a \ast b)\rightarrow
c$,
\item[{\rm (I4)}] $({\mathop\bigvee\limits _{j\in J}a _j)}\rightarrow b={\mathop\bigwedge\limits _{j\in J}}(a _j\rightarrow
b)$,

\item[{\rm (I5)}] $a \rightarrow({\mathop\bigwedge\limits _{j\in J}b_j)}={\mathop\bigwedge\limits _{j\in J}}(a \rightarrow
b_j)$.

\end{enumerate}

$L=(L,\ast)$ is called  {\it a frame or a complete Heyting algebra } when $\ast=\wedge$.

$L=(L,\ast)$ is called {\it a complete MV-algebra} if it holds that {\bf (MV)}\ \ $\forall a,b \in L, a\vee b=(a\rightarrow b)\rightarrow b$.

%
%
%
%
%
%
%

\subsection{$L$-fuzzy set over group}

By an $L$-fuzzy set in $X$ we mean a mapping $\lambda:X\rightarrow L$, the set of all $L$-fuzzy sets in $X$ is denoted by $L^X$. The operations
$\vee,\wedge,\ast,\rightarrow$ on $L$ can be defined onto $L^X$ pointwisely.

Let  $f:X\longrightarrow Y$ be a mapping. Then the mappings
$f^{\rightarrow}: L^X\longrightarrow L^Y$ and $f^{\leftarrow}:
L^Y\longrightarrow L^X$ determined by
$$\forall \lambda\in L^X, y\in Y, f^{\rightarrow}(\lambda)(y)=\bigvee_{f(x)=y}\lambda(x);\ \ \forall \nu\in
L^Y, x\in X, f^{\leftarrow}(\nu)(x)=\nu (f(x))$$ are called the Zadeh mappings.

For a group $(X,\cdot_X)$,  we denote the group operations by $$m_X:X\times X\longrightarrow
X, (x,y)\mapsto x\cdot_X y\ (xy\ {\rm for\ simplicity}); \ r_X:X\longrightarrow X, x\mapsto x^{-1},$$ and {\bf if there is no confusion we will short $\cdot_X$, $m_X$ and $r_X$ as $\cdot$, $m$ and $r$.}

%

A mapping  $f:(X,\cdot)\longrightarrow (Y,\cdot)$ is called a group
homomorphism if $f(xz)=f(x)f(z)$ for all $x,z\in X$.

 For any $\lambda,\mu\in L^X$,  define
$\lambda^{-1}, \lambda\odot\mu\in L^X, \lambda\times \mu\in L^{X\times X}$ by
$$\forall z\in X, \lambda^{-1}(z)=\lambda(z^{-1}), \ \ (\lambda\odot\mu)(z)=\bigvee_{xy=z}\big(\lambda(x)\wedge \mu(y)\big), \forall x=(x_1,x_2)\in X\times X, \ \  (\lambda\times
\mu)(x)=\lambda(x_1)\wedge \mu(x_2).$$

For $x\in X$, we also use $x$ to denote its characteristic mapping $\top_{\{x\}}$ for simplicity.

%

For $\lambda,\mu\in L^X$, we define $S_X(\lambda,\mu)= {\mathop \bigwedge\limits_{x\in
X}}(\lambda(x)\rightarrow \mu(x))$. This induces a mapping $S_X:L^X\times L^X\longrightarrow L^X$, called the subsethood degree mapping
\cite{R.B02,DZ07}.

The following two lemmas collect the fundamental properties of the considered operations.

\begin{lemma} (\cite{TG142,TA17}) \label{lemma01} Let $e$ be the identity element
of a group $(X,\cdot)$  and $x,y\in X$,  $\lambda,\mu \in L^X$.

{\rm (1)}  $\forall z\in X,x\odot\mu(z)=\mu(x^{-1}z)$ and  $e\odot \lambda=\lambda$,




%
 {\rm (2)} $(x\odot
\lambda)\odot \mu=x\odot (\lambda\odot \mu)$,


{\rm (3)} $(xy)\odot \lambda=(x\odot y) \odot \lambda=x\odot (y \odot \lambda)$,

%
{\rm (4)} $\forall \lambda_i,\mu_i\in L^X, i=1,2$,  $(\lambda_1\wedge\lambda_2)\odot (\mu_1\wedge\mu_2)\leq (\lambda_1\odot
\mu_1)\wedge(\lambda_2\odot\mu_2)$.
\end{lemma}

\begin{lemma} (\cite{R.B02,fang2017}) \label{lemmaprd} Let  $\lambda_i,\mu_i (i=1,2) \subseteq L^X$.





{\rm (1)} $S_{X}(\lambda_1,\mu_1)\wedge S_{X}(\lambda_2,\mu_2)\leq
S_{X\times X}(\lambda_1\times \lambda_2, \mu_1\times \mu_2)$,

 {\rm (2)} $\forall x\in X, S_X(\lambda_1,\mu_1)\leq S_X(x\odot\lambda_1, x\odot\mu_1)$,

{\rm (3)} $S_{X}(\lambda_1,\mu_1)\wedge
S_{X}(\lambda_2,\mu_2)\leq S_{X}(\lambda_1\odot \lambda_2, \mu_1\odot \mu_2)$.

\end{lemma}


%
%
%
%
%

\subsection{$\top$-filters}

In this subsection, we summarize some necessary properties about $\top$-filters.

%

\begin{definition} (\cite{U.H99}) A nonempty subset $\mathbb{F}\subseteq L^X$ is called a $\top$-filter on  $X$ provided that:

(TF1) $\forall \lambda\in
\mathbb{F}$, $\bigvee\limits_{x\in X}\lambda(x)=\top$,

(TF2) $\forall \lambda,\mu\in
\mathbb{F}$, $\lambda\wedge \mu\in \mathbb{F}$,

(TF3) $\forall \lambda\in L^X$,  $\bigvee\limits_{\mu\in
\mathbb{F}}S_X(\mu,\lambda)=\top\Longrightarrow\lambda\in\mathbb{F}$.

The set of all $\top$-filters on $X$ is denoted by
$\mathbb{F}_L^{\top}(X)$.
 $\top$-filter is an important extension of the classical filter.

\end{definition}

\begin{definition} (\cite{U.H99}) A nonempty subset $\mathbb{B}\subseteq L^X$ is called a $\top$-filter base on $X$ provided that:

(TB1) $\forall \lambda\in
\mathbb{B}$, $\bigvee\limits_{x\in X}\lambda(x)=\top$.

(TB2)  $\forall \lambda,\mu\in \mathbb{B}$, $\bigvee\limits_{\nu\in
\mathbb{B}}S_X(\nu,\lambda\wedge\mu)=\top$.
\end{definition}

Each $\top$-filter base $\mathbb{B}$ generates a $\top$-filter $\mathbb{F}_{\mathbb{B}}:=\{\lambda\in L^X|\bigvee\limits_{\mu\in \mathbb{B}}S_X(\mu,\lambda)=\top\}.$

\begin{example} (\cite{U.H99,LZ23}) \label{lemimge} Let $f:X\longrightarrow Y$ be a mapping.

(1) For any $\mathbb{F}\in \mathbb{F}_L^{\top}(X)$, the $\top$-filter generated by $\{f^\rightarrow(\lambda)|\lambda \in \mathbb{F}\}\subseteq L^Y$
as a $\top$-filter base, is denoted by
$f^\Rightarrow(\mathbb{F})$.  Obviously, $\mu\in f^\Rightarrow(\mathbb{F})\Longleftrightarrow f^\leftarrow(\mu)\in \TF$. Furthermore, if $\TB$ is a $\top$-filter base of $\TF$, then the family $\{f^\rightarrow(\lambda)|\lambda \in \mathbb{B}\}$  forms a $\top$-filter base of $f^\Rightarrow(\TF)$.

(2) For any $\mathbb{G}\in \mathbb{F}_L^{\top}(Y)$,  $\{f^\leftarrow(\mu)|\mu \in \mathbb{G}\}\subseteq L^X$ forms a
$\top$-filter base $\Longleftrightarrow$ $\forall \mu\in \mathbb{G}$, $\bigvee\limits_{y\in
f(X)}\mu(y)=\top$.  The $\top$-filter generated by it (if exists) is denoted as $f^\Leftarrow(\mathbb{G})$. Note that $f^\Leftarrow(\mathbb{G})$ always exist provided that $f$ is surjective. In addition, if $\TB$ is a $\top$-filter base of $\TG$ and $f^\Leftarrow(\mathbb{G})$ exists then the family $\{f^\leftarrow(\lambda)|\lambda \in \mathbb{B}\}$  forms a $\top$-filter base of $f^\Leftarrow(\mathbb{G})$.

{\rm (3)} For any $x\in X$, the family $[x]_{\top}:=\{\mu\in L^X|\mu(x)=\top\}$ is a $\top$-filter on $X$ with a $\top$-filter base $\TB=\{\top_{\{x\}}\}$.

\end{example}

\begin{definition} (\cite{Yu17}) \label{defntfb} Let $\mathbb{F}_i\in \mathbb{F}_L^{\top}(X_i)$$(i=1,2)$. Then the family $\{\lambda_1\times \lambda_2|\lambda_i\in \mathbb{F}_i\}$ forms a $\top$-filter base on $X_1\times X_2$, and the generated  $\top$-filter denoted by $\mathbb{F}_1\times \mathbb{F}_2$. Furthermore,  if  $\mathbb{B}_i$ is a $\top$-filter base of $\mathbb{F}_i$ then  $\mathbb{B}=\{\lambda_1\times \lambda_2|\lambda_i\in \mathbb{B}_i\}$ is a  $\top$-filter base of $\mathbb{F}_1\times \mathbb{F}_2$.\end{definition}

Let $\mathbb{F}_i\in \mathbb{F}_L^{\top}(X_i)$$(i=1,2)$ and $p_i:X_1\times X_2\longrightarrow X_i$ be the projections. It is easily seen that $\forall \lambda_i\in L^{X_i}$,  $p_1^{\leftarrow}(\lambda_1)\wedge p_2^{\leftarrow}(\lambda_2)=\lambda_1\times \lambda_2$. Hence  $\mathbb{F}_1\times \mathbb{F}_2$ has the $\top$-filter base $\TB=\{p_1^{\leftarrow}(\lambda_1)\wedge p_2^{\leftarrow}(\lambda_2)|\lambda_i\in \TF_i\}.$

\begin{proposition} (\cite{YY22}) \label{TPD} Let $\TF,\TG\in \TF^\top_L(X)$ and $\mathbb{K}\in \TF_L^\top(X\times X)$.

{\rm (1)}  $p_1^\Rightarrow(\TF\times \TG)=\TF$ and $p_2^\Rightarrow(\TF\times \TG)=\TG$,

{\rm (2)}  $\mathbb{K}\supseteq p_1^\Rightarrow(\mathbb{K})\times p_2^\Rightarrow(\mathbb{K})$.
\end{proposition}

%


\subsection{Categoric theory}

For general categoric theory, please refer to \cite{AHS}. We only recall the notion of Cartesian closedness.

A category {\bf A} is called Cartesian closed if it satisfies the following conditions:

(1) For a pair of {\bf A}-objects $(X,Y)$ there exists  a product $X\times Y$ in {\bf A}.

(2) For a pair of {\bf A}-objects $(X,Y)$, there is some {\bf A}-object $Y^X$ (called power
object) and some {\bf A}-morphism $ev : Y^ X \times  X \longrightarrow Y$ (called evaluation
morphism) s.t. for any {\bf A}-object $Z$ and any {\bf A}-morphism $\psi: Z\times X\longrightarrow Y$, there is a unique {\bf A}-morphism $\psi^\lozenge: Z \longrightarrow Y^X$ s.t. $ev\circ (\psi^\lozenge \times id_X) = \psi$.

\subsection{$\top$-convergence spaces and $\top$-convergence groups}

In the subsection, we review some outcomes of  $\top$-convergence spaces and $\top$-convergence groups.

\begin{definition} (\cite{fang2017}) A $\top$-convergence space is a pair $(X, \mathbb{C}_X)$ (when there is no confusion, we omit the subscript $_X$ and simplify $\mathbb{C}_X$ as $\TC$), where $\TC: \TF_L(X)\longrightarrow 2^X$ is a mapping satisfying:

(TC1) $\forall x\in X$, $[x]_{\top}\stackrel{\TC}{\dashrightarrow} x$,

(TC2) $\forall x\in X$, $\TF\subseteq \TG, \TF\stackrel{\TC}{\dashrightarrow} x\Longrightarrow \TG\stackrel{\TC}{\dashrightarrow} x$,

\noindent where $\TF\stackrel{\TC}{\dashrightarrow} x$ means that $x\in \TC(\TF)$. The pair $(X,\TC)$ is called a $\top$-convergence space ($\top$-Cons for short).
\end{definition}

A mapping  $(X,\TC_X)\stackrel{f}\longrightarrow (Y,\TC_Y)$ between $\top$-Cons is said to be continuous at if  $\forall x\in X,  \forall \TF\in \TF_L^\top(X)$, $\TF\stackrel{\TC_X}{\dashrightarrow} x$ ensures $f^\Rightarrow(\TF)\stackrel{\TC_Y}{\dashrightarrow} f(x)$. Let {\bf TCons} denote the category of $\top$-Cons and
continuous mappings.

A source $(X\stackrel{f_j}{\longrightarrow} (X_j,\TC_j))_{j\in J}$ in {\bf TCons} consists of a family of mappings $(X\stackrel{f_j}{\longrightarrow} X_j)_{j\in J}$ together with a family of $\top$-Cons $(X_j,\TC_j)_{j\in J}$. The {\it initial structure} w.r.t. the source $\big(X\stackrel{f_j}{\longrightarrow} (X_j,\TC_j)\big)_{j\in J}$ is meant the $\top$-convergence structure $\TC$ on $X$ satisfying  the following initial conditions:

(1) Each $f_j:(X,\TC)\longrightarrow (X_j,\TC_j)$ ($j\in J$) is continuous,

(2)  $\forall (Y,\TV)\in ${\bf TCons},  a mapping $(Y,\TV)\stackrel{ f}\longrightarrow(X,\TC)$ is continuous iff all $(Y,\TV)\stackrel{f_j\circ f}\longrightarrow (X_j,\TC_j)$ are continuous.

\begin{theorem} (\cite{Yu17}) \label{tinitial} Each source $\big(X\stackrel{f_j}{\longrightarrow} (X_j,\TC_j)\big)_{j\in J}$ in {\bf TCons} has an initial structure $\TC$ defined by $\forall x\in X, \forall \TF\in \TF_L^\top(X)$, $\TF\stackrel{\TC}{\dashrightarrow} x$ iff $\forall j\in J, f^\Rightarrow_j(\TF)\stackrel{\TC_j}{\dashrightarrow} f_j(x)$. Hence {\bf TCons} is a topological category.
\end{theorem}


%
%

Because every topological category has  products, the category {\bf TCons} has products, and therefore has finite products.

\begin{definition} (\cite{Yu17}) Let $(X,\TC)$ be a $\top$-Cons and $\TC\times
\TC$ be the initial structure of the source $\big(X\times X\stackrel{p_i}\longrightarrow (X,\TC)\big)_{i=1,2}$ in {\bf TCons}. Then we call  the pair $(X\times X,\TC\times\TC)$  the product space of $(X,\TC)$. In precise, for any $\TF\in \TF_L^\top(X\times X)$ and any $x=(x_1,x_2)\in X\times X$, $\TF\stackrel{\TC\times \TC}{\dashrightarrow} x\Longleftrightarrow p_i(\TF)\stackrel{\TC}{\dashrightarrow} x_i, i=1,2.$
\end{definition}

\begin{proposition} \label{conpr} Let $f_i:(X_i, C_{X_i})\longrightarrow (Y_i, C_{Y_i})$ ($i=1,2$) be continuous mappings between $\top$-convergence spaces $(X_i, C_{X_i})$ and $(Y_i, C_{Y_i})$. Define $f_1\times f_2$: $X_1\times X_2\longrightarrow Y_1\times Y_2$ by $\forall (x_1,x_2)\in X_1\times X_2$, $f_1\times f_2 (x_1,x_2)=(f_1(x_1),f_2(x_2))$, then $f_1\times f_2$ is a continuous mapping between the corresponding product spaces.
\end{proposition}

\begin{proof} Notice that each $f_i\circ p_i$ is continuous. Then by $p_i\circ (f_1\times f_2)=f_i\circ p_i$ we obtain that each $p_i\circ (f_1\times f_2)$ is continuous. It follows by the initial condition (2) of product space we have $f_1\times f_2$ is continuous.
\end{proof}

\begin{definition} (\cite{fang2017}) (1) A $\top$-Cons $(X,\TC)$ is called a $\top$-limit space whenever it fulfills

(LT)  $\TF, \TG\stackrel{\TC}{\dashrightarrow} x\Longrightarrow \TF\cap \TG\stackrel{\TC}{\dashrightarrow} x$.

(2) A $\top$-Cons $(X,\TC)$ is called pretopological whenever it fulfills

(PT)  $\forall x\in X$, $\mathbb{U}_x\stackrel{\TC}{\dashrightarrow} x$, where $\mathbb{U}_x=\bigcap\{\TF:\TF \stackrel{\TC}{\dashrightarrow} x\}$.

(3) A pretopological $\top$-Cons $(X,\TC)$ is called topological whenever it fulfills

(TT)  $\forall x\in X, \forall \lambda\in \mathbb{U}_x$, there exists $\lambda^\ast\in \mathbb{U}_x$ with $\lambda^\ast\leq \lambda$ and for each $z\in X$, there exists $\lambda_z\in \mathbb{U}_z$ s.t. $\lambda^\ast(z)\leq  S_X(\lambda_z, \lambda)$.
\end{definition}

\begin{definition} (\cite{LZ23}) Let $(X,\TC)$
be a  $\top$-Cons over a group $(X,\cdot)$. Then we call the triple
$(X,\cdot,\TC)$  a  $\top$-convergence
group  whenever  the group operations $(X\times X, \TC\times
\TC)\stackrel{m}{\longrightarrow} (X,\TC)$ and $(X,\TC)\stackrel{r}{\longrightarrow}
(X,\TC)$ are continuous.
\end{definition}

In the following, we use  {\bf TConG} to denote the category of  $\top$-convergence groups and
continuous group homomorphisms.

\section{The characterization via $\odot$-product of $\top$-filters and the localization}

In this section, we will present a new characterization on $\top$-convergence groups by the $\odot$-product of $\top$-filters, and then prove that $\top$-convergence groups can be
localizable.

At first, we fix the notion of $\odot$-product of $\top$-filters.

\begin{proposition} \label{proppc} Let $\mathbb{F}_i\in \mathbb{F}_L^{\top}(X)$$(i=1,2)$. Then the family $\{\lambda\odot\mu|\lambda\in \TF_1, \mu\in \TF_2\}$ forms a $\top$-filter base on $X$, the resulted $\top$-filter is denoted as $\TF_1\odot\TF_2$. Furthermore, if $\mathbb{B}_i$ is a $\top$-filter base of $\mathbb{F}_i$ then the family  $\mathbb{B}=\{\lambda_1\odot \lambda_2|\lambda_i\in \mathbb{B}_i\}$ is a  $\top$-filter base of $\mathbb{F}_1\odot \mathbb{F}_2$.
\end{proposition}

\begin{proof} We only check the first part and that the second part is obviously.

For any $\lambda\in \TF_1, \mu\in \TF_2$, it holds that $\bigvee\limits_{x\in X}\lambda(x)=\top$ and $\bigvee\limits_{y\in X}\mu(y)=\top$. Then
$$\bigvee\limits_{z\in X}(\lambda\odot \mu)(z)=\bigvee\limits_{z\in X}\bigvee\limits_{xy=z}(\lambda(x)\wedge \mu (y))=\bigvee\limits_{x,y\in X}(\lambda(x)\wedge \mu (y))\geq \bigvee\limits_{x,y\in X}(\lambda(x)\ast \mu (y))=\Big(\bigvee\limits_{x\in X}\lambda(x)\Big) \ast \Big(\bigvee\limits_{y\in X}\mu(y)\Big)=\top.$$Hence (TB1) is satisfied.

 For any $\lambda_1,\lambda_2\in \TF_1, \mu_1,\mu_2\in \TF_2$, it follows by   Lemma \ref{lemma01} (4) that $(\lambda_1\wedge\lambda_2)\odot (\mu_1\wedge\mu_2)\leq (\lambda_1\odot
\mu_1)\wedge(\lambda_2\odot\mu_2), $ and then by  $\lambda_1\wedge\lambda_2\in \TF_1, \mu_1\wedge\mu_2\in \TF_2$ we get
$$\bigvee\limits_{\lambda_3\in \TF_1, \mu_3\in \TF_2}S_X\Big(\lambda_3\odot \mu_3, (\lambda_1\odot
\mu_1)\wedge(\lambda_2\odot\mu_2)\Big)\geq S_X\Big((\lambda_1\wedge\lambda_2)\odot (\mu_1\wedge\mu_2), (\lambda_1\odot
\mu_1)\wedge(\lambda_2\odot\mu_2)\Big)=\top.$$Hence (TB2) is satisfied.
\end{proof}

\begin{definition} (1) Let $\mathbb{F}, \TG\in \mathbb{F}_L^{\top}(X)$. Then the $\top$-filter $\TF\odot\TG$ is called the $\odot$-product of $\TF$ and $\TG$.

(2) Let $\mathbb{F}\in \mathbb{F}_L^{\top}(X)$. Then it is easily proved that the family $\{\lambda^{-1}|\lambda\in \TF\}$ forms a $\top$-filter on $X$, denoted as $\TF^{-1}$, called the inverse of $\TF$.
\end{definition}

We collect some properties of $\odot$-product of $\top$-filters for later use.

\begin{proposition} \label{TFP} Let $\TF,\TG, \mathbb{H}, \in \TF^\top_L(X)$. 

{\rm (1)}  $m^\Rightarrow(\TF\times \TG)=\TF\odot \TG$,

{\rm (2)}  $[e]_\top\odot \TF=\TF= \TF\odot[e]_\top$,

{\rm (3)}  $(\TF\odot \TG)\odot\mathbb{H}=\TF\odot (\TG\odot\mathbb{H})$,

{\rm (4)}  $\forall x,y\in X$, $[x]_\top\odot [y]_\top= [xy]_\top$.

\end{proposition}

\begin{proof} (1) From Example \ref{lemimge} (1) and Lemma \ref{lemma01} (3) we know that $m^\Rightarrow(\TF\times \TG)$ and $\TF\odot \TG$ have a common $\top$-filter base  $\{m^\rightarrow(\lambda\times \mu)=\lambda\odot\mu|\lambda\in \TF, \mu\in \TG\}$, and hence they are equal.

(2) We only prove $[e]_\top\odot \TF=\TF$. The other equality is similar. For any $\lambda \in \TF$, from $\lambda=e\odot \lambda$ we conclude that $\TF\subseteq [e]_\top\odot \TF$. Conversely, for any $\mu\in [e]_\top, \lambda \in \TF$, from $\lambda=e\odot\lambda\leq \mu\odot \lambda$ we get that $\mu\odot \lambda\in \TF$, and so $[e]_\top\odot \TF\subseteq\TF$.

(3) It follows by that $\{\lambda\odot \mu\odot \nu|\lambda\in \TF, \mu\in \TG,\nu\in \mathbb{H}\}$ forms a common $\top$-filter base of $(\TF\odot \TG)\odot\mathbb{H}$ and $\TF\odot (\TG\odot\mathbb{H})$.

(4) For any $\lambda\in [x]_\top, \mu\in [y]_\top$, we have $(\lambda\odot \mu)(xy)\geq \lambda(x)\wedge \mu(y)=\top$, which means $\lambda\odot \mu\in [xy]_\top$, and so  $[x]_\top\odot [y]_\top\subseteq [xy]_\top$. Conversely, let $\lambda \in [xy]_\top$, then by $\top_{\{x\}}\in [x]_\top$, $\top_{\{y\}}\in [y]_\top$ and $\lambda\geq \top_{\{x\}}\odot \top_{\{y\}}$ we obtain $\lambda\in [x]_\top\odot [y]_\top$, and so $ [xy]_\top\subseteq [x]_\top\odot [y]_\top$.\end{proof}
%

The theorem below presents an important characterization on $\top$-convergence group.

\begin{theorem}\label{thmct2} (Characterization via $\odot$-Product of $\top$-filters) Let $(X,\TC)$
be a  $\top$-Cons over a group $(X,\cdot)$. Then $(X,\cdot,\TC)$ is a $\top$-convergence group iff the following two conditions are satisfied:

{\rm (TCG1)} $\forall \TF,\TG\in \TF_L^\top(X)$, $\forall x,y\in X$, $\TF\stackrel{\TC}{\dashrightarrow} x, \TG\stackrel{\TC}{\dashrightarrow} y$ implies $\TF\odot \TG\stackrel{\TC}{\dashrightarrow} xy$,

{\rm (TCG2)} $\forall \TF\in \TF_L^\top(X)$, $\forall x\in X$, $\TF\stackrel{\TC}{\dashrightarrow} x$ implies $\TF^{-1}\stackrel{\TC}{\dashrightarrow} x^{-1}$.
\end{theorem}

\begin{proof} Obviously, $r$ is continuous iff (TCG2) holds. Next, we prove that $m$ is continuous iff (TCG1) holds.

Assume that $m$ is continuous. Take any $\TF\stackrel{\TC}{\dashrightarrow} x, \TG\stackrel{\TC}{\dashrightarrow} y$, then $\TF\times \TG \stackrel{\TC\times \TC}{\dashrightarrow} (x,y)$ by the definition of product space and Proposition \ref{TPD} (1). It follows by the continuity of $m$ and Proposition \ref{TFP} (1) that  $\TF\odot \TG=m^\Rightarrow(\TF\times \TG)\stackrel{\TC}{\dashrightarrow} xy$. Hence (TCG1) holds.

Conversely, let (TCG1) hold.  Take any $\TF \stackrel{\TC\times \TC}{\dashrightarrow} x=(x_1,x_2)$, then $p_i^\Rightarrow(\TF)\stackrel{\TC}{\dashrightarrow} x_i$ ($i=1,2$). It follows by  Proposition \ref{TPD} (2),  Proposition \ref{TFP} (1) and (TCG1) that $$\TF\supseteq p_1^\Rightarrow(\TF)\times p_2^\Rightarrow(\TF),\ {\rm and}\ m^\Rightarrow( p_1^\Rightarrow(\TF)\times p_2^\Rightarrow(\TF))= p_1^\Rightarrow(\TF)\odot p_2^\Rightarrow(\TF) \stackrel{\TC}{\dashrightarrow} x_1x_2,$$ which means $m^\Rightarrow(\TF)\stackrel{\TC}{\dashrightarrow} m(x)$. Hence $m$ is continuous.
\end{proof}

%



%
%
%


\begin{definition} A $\top$-convergence group $(X,\cdot, \TC)$ is called localizable if the associated $\top$-convergence structure $\TC$ is determined by the convergence at the identity $e$.
\end{definition}

The following theorem shows that a $\top$-convergence group is localizable.

\begin{theorem} (Localizable theorem) \label{theorem01} For $(X,\cdot, \TC)\in $ {\bf TConG}, then$$\TF\stackrel{\TC}{\dashrightarrow} x\Longleftrightarrow [x^{-1}]_\top\odot \TF\stackrel{\TC}{\dashrightarrow} e\Longleftrightarrow \TF\odot[x^{-1}]_\top\stackrel{\TC}{\dashrightarrow} e.$$
\end{theorem}
\begin{proof} We prove $\TF\stackrel{\TC}{\dashrightarrow} x\Longleftrightarrow [x^{-1}]_\top\odot \TF\stackrel{\TC}{\dashrightarrow} e$ as example, the other equivalence is similar.

$\Longrightarrow$. From $\TF\stackrel{\TC}{\dashrightarrow} x$, $[x^{-1}]_\top\stackrel{\TC}{\dashrightarrow} x^{-1}$ and (TCG1) we obtain that $[x^{-1}]_\top\odot\TF \stackrel{\TC}{\dashrightarrow} x^{-1}x=e$.

$\Longleftarrow$. From $[x^{-1}]_\top\odot\TF \stackrel{\TC}{\dashrightarrow} e$,  $[x]_\top\stackrel{\TC}{\dashrightarrow} x$, (TCG1) and Proposition \ref{TFP} (2)-(4) we obtain that  $\TF=[x]_\top \odot[x^{-1}]_\top\odot\TF\stackrel{\TC}{\dashrightarrow} xe=x$.
\end{proof}

\section{Pretopological $\top$-convergence group is equivalent to topological $\top$-convergence group}

A pretopological $\top$-convergence space is typically not a topological $\top$-convergence space. In this section, we will introduce a category of pretopological $\top$-convergence groups, and demonstrate that the $\top$-convergence space associated with a pretopological $\top$-convergence group is, in fact, a topological $\top$-convergence space. Consequently, this establishes that a pretopological $\top$-convergence group is equivalent to a topological $\top$-convergence group.

\begin{definition} (\cite{LZ23}) A  $\top$-convergence group $(X,\cdot,\TC)$ is called a topological $\top$-convergence  group (resp., $\top$-limit group) provided that the $\top$-convergence space $(X,\TC)$ is a topological $\top$-convergence space (resp., a $\top$-limit space).
\end{definition}

Naturally, between $\top$-limit groups and topological $ \top $-convergence groups, we can define a new subcategory of $\top$-convergence groups as below.

\begin{definition} A  $\top$-convergence group $(X,\cdot,\TC)$ is called a pretopological $\top$-convergence  group provided  that $(X,\TC)$ is a  pretopological $\top$-convergence space.
\end{definition}

\begin{example} For a group $(X,\cdot)$, let $\TC_{ds}$ (resp., $\TC_{ids}$) be defined by $\forall x\in X, \forall \TF\in \TF_L^\top(X)$,  $$\TF\stackrel{\TC_{ds}}{\dashrightarrow} x\Longleftrightarrow \TF\supseteq [x]_{\top}\ \  {\rm (resp.,}\ \TF\stackrel{\TC_{ids}}{\dashrightarrow} x\Longleftrightarrow \TF\in \TF_L^\top(X)). $$It is easy to prove that $(X,\cdot, \TC_{ds})$ (resp., $(X,\cdot, \TC_{ids})$) constitutes  a pretopological  $\top$-convergence group.
\end{example}

The next proposition collects the properties of the $\top$-filter $\TU_x$ in a pretopological $\top$-convergence group.

\begin{proposition} \label{propun} Let $(X,\cdot, \TC)$ be a pretopological $\top$-convergence group.

(1) $\lambda\in \TU_x$ iff $x^{-1}\odot \lambda \in \TU_e$.

(2) $\lambda\in \TU_e$ iff $x\odot \lambda \in \TU_x$.

(3) $\TU_e\subseteq \TU_e\odot \TU_e$.
\end{proposition}

\begin{proof} (1) $\Longrightarrow$. Let $x\in X$, recall that $\TU_x=\bigcap\{\TF:  \TF\stackrel{\TC}{\dashrightarrow} x\}$. For any $\TG\stackrel{\TC}{\dashrightarrow} e$, one get $[x^{-1}]_\top\odot ([x]_\top\odot \TG)=\TG\stackrel{\TC}{\dashrightarrow} e$,  then from the Localizable theorem (i.e., Theorem \ref{theorem01}) we obtain $[x]_\top\odot \TG \stackrel{\TC}{\dashrightarrow} x$, so $\TU_x\subseteq [x]_\top\odot \TG$ by the pretopological condition (PT). So for any $\lambda\in \TU_x$, one get $\lambda\in [x]_\top\odot \TG$, hence by Proposition \ref{proppc}
\begin{eqnarray*}\top&=&\bigvee_{\nu\in \TG}S_X(x\odot \nu, \lambda)\leq\bigvee_{\nu\in \TG}S_X(x^{-1}\odot x\odot \nu, x^{-1}\odot\lambda)=\bigvee_{\nu\in \TG}S_X(\nu, x^{-1}\odot\lambda)
\end{eqnarray*} which means $x^{-1}\odot\lambda\in \TG$ for any $\TG\stackrel{\TC}{\dashrightarrow} e$, so $x^{-1}\odot\lambda\in \TU_e$, as desired.

$\Longleftarrow$. Let $x^{-1}\odot\lambda\in \TU_e$ and $\TF\stackrel{\TC}{\dashrightarrow} x$.  Then from Theorem \ref{theorem01} we obtain $[x^{-1}]_\top\odot \TF\stackrel{\TC}{\dashrightarrow} e$, so $x^{-1}\odot \lambda \in \TU_e\subseteq [x^{-1}]_\top\odot \TF$ by the pretopological condition (PT). Hence $\lambda=x\odot x^{-1}\odot \lambda \in  [x]_\top\odot [x^{-1}]_\top\odot \TF=\TF$ for any $\TF\stackrel{\TC}{\dashrightarrow} x$. So, $\lambda\in \TU_x$ as desired.

(2) It follows by (1).

(3) From the pretopological condition (PT) on have $\TU_e\stackrel{\TC}{\dashrightarrow} e$, then  $\TU_e\odot \TU_e\stackrel{\TC}{\dashrightarrow} e$ by Theorem \ref{thmct2}. Hence $\TU_e\subseteq \TU_e\odot \TU_e$.
\end{proof}

 Let $(X,\cdot, \TC)$ be a pretopological $\top$-convergence group. For any $\lambda\in L^X$, we define $\lambda^\ast\in L^X$ by $\forall x\in X$, $$\lambda^\ast(x)=\bigvee_{\mu\in \TU_x}S_X(\mu, \lambda).$$

The theorem below shows that the pretopological $\top$-convergence groups are equivalent to the topological $\top$-convergence groups. Hence, the topological condition (TT) in the definition of topological $\top$-convergence space can be removed.

\begin{theorem} \label{themre} Let $(X,\cdot, \TC)$ be a pretopological $\top$-convergence group and $L$ be an  MV-algebra. Then  $(X, \TC)$ is a topological $\top$-convergence space.
\end{theorem}

\begin{proof} We prove that $(X, \TC)$ satisfies the topological condition (TT). We will complete the proof in the following three steps.

(1) For any $\lambda\in\mathbb{U}_x$, we prove that $\lambda^\ast\in\mathbb{U}_x$.

From Proposition  \ref{propun}  we have $x^{-1}\odot \lambda\in \mathbb{U}_e\odot \mathbb{U}_e$, so $\bigvee_{\lambda_1,\lambda_2\in \mathbb{U}_e} S_X(\lambda_1\odot\lambda_2, x^{-1}\odot
\lambda)=\top.$ Notice that
\begin{eqnarray*}S_X\Big(\lambda_1\odot\lambda_2, x^{-1}\odot \lambda\Big)&\leq&S_X\Big((x\odot\lambda_1)\odot\lambda_2, \lambda\Big)\\
&=&\bigwedge_{y\in
X}\Big(\big[(x\odot\lambda_1)\odot\lambda_2\big](y)\rightarrow \lambda(y)\Big)\\
&=&\bigwedge_{y\in X}\Big(\bigvee_{z\in
X}\big[(x\odot\lambda_1)(z)\wedge\lambda_2(z^{-1}y)\big]\rightarrow \lambda(y)\Big)\\
&\leq&\bigwedge_{y\in X}\Big(\bigvee_{z\in
X}\big[(x\odot\lambda_1)(z)\ast\lambda_2(z^{-1}y)\big]\rightarrow \lambda(y)\Big)\\
&=&\bigwedge_{y\in X}\bigwedge_{z\in
X}\Big(\big[(x\odot\lambda_1)(z)\ast\lambda_2(z^{-1}y)\big]\rightarrow \lambda(y)\Big)\\
&=&\bigwedge_{y,z\in X}\Big(\big[x\odot\lambda_1\big](z)\rightarrow\big[\lambda_2(z^{-1}y)\rightarrow
\lambda(y)\big]\Big)\\
&=&\bigwedge_{y,z\in X}\Big(\big[x\odot\lambda_1\big](z)\rightarrow\big[(z\odot\lambda_2)(y)\rightarrow
\lambda(y)\big]\Big)\\
&=&\bigwedge_{z\in X}\Big(\big[x\odot\lambda_1\big](z)\rightarrow\bigwedge_{y\in
X}\big[(z\odot\lambda_2)(y)\rightarrow \lambda(y)\big]\Big)\\
&=&\bigwedge_{z\in X}\Big(\big[x\odot\lambda_1\big](z)\rightarrow
S_X(z\odot\lambda_2, \lambda)\Big) \ \ {\rm by}\ x\odot\lambda_1\in \mathbb{U}_x, z\odot\lambda_2\in \mathbb{U}_z\\
&\leq& \bigvee_{\nu\in \mathbb{U}_x}\bigwedge_{z\in X}\Big(\nu(z)\rightarrow
\bigvee_{\mu\in \mathbb{U}_z} S_X(\mu,
\lambda)\Big)\\
&=& \bigvee_{\nu\in \mathbb{U}_x}\bigwedge_{z\in X}\Big(\nu(z)\rightarrow \lambda^\ast (z)\Big)
\\
&=& \bigvee_{\nu\in \mathbb{U}_x}S_X\Big(\nu, \lambda^\ast\Big).
\end{eqnarray*}

Hence
\begin{eqnarray*}\top&=&\bigvee_{\lambda_1,\lambda_2\in \mathbb{U}_e} S_X(\lambda_1\odot\lambda_2, x^{-1}\odot
\lambda)\leq \bigvee_{\nu\in \mathbb{U}_x}S_X(\nu,\lambda^\ast),
\end{eqnarray*}
which means $\lambda^\ast\in \mathbb{U}_x$.

(2) For any $\lambda\in\mathbb{U}_x$, we prove that $\lambda^\ast \leq \lambda$.

For any $y\in X, \mu\in \TU_y$, from $\TU_y\subseteq [y]_\top$, we get $\mu(y)=\top$. Hence $$\lambda^\ast(y)=\bigvee_{\mu\in \TU_y}S_X(\mu, \lambda)\leq \bigvee_{\mu\in \TU_y}(\mu(y)\rightarrow \lambda(y))=\lambda(y).$$

(3)  For any $\lambda\in\mathbb{U}_x$ and for each $y\in X$, we define $\lambda_y\in L^X$ by $\lambda_y(z)=\lambda^\ast(y)\rightarrow \lambda(z)$. Then we prove that $\lambda_y\in \mathbb{U}_y$ and $\lambda^\ast(y)\leq  S_X(\lambda_y, \lambda)$.
$$S_X(\lambda_y, \lambda)=\bigwedge_{z\in X}(\lambda_y(z)\rightarrow \lambda(z))=\bigwedge_{z\in X}((\lambda^\ast(y)\rightarrow \lambda(z))\rightarrow \lambda(z))\geq \lambda^\ast(y).$$
\begin{eqnarray*}\bigvee_{\mu\in \TU_y}S_X(\mu, \lambda_y)&=&\bigvee_{\mu\in \TU_y}\bigwedge_{z\in X}(\mu(z)\rightarrow (\lambda^\ast(y)\rightarrow \lambda(z)))\\
&=&\bigvee_{\mu\in \TU_y}(\lambda^\ast(y) \rightarrow S_X(\mu, \lambda))\ \ {\rm by\ MV}\\
&=&\lambda^\ast(y) \rightarrow \bigvee_{\mu\in \TU_y}(S_X(\mu, \lambda))\\
&=&\lambda^\ast(y) \rightarrow \lambda^\ast(y)=\top,
\end{eqnarray*}so $\lambda_y\in \TU_y$.

(1)-(3) shows  that $(X, \TC)$ satisfies the topological condition (TT).
\end{proof}

\begin{remark} In \cite{LZ23}, Zhang and Pang verified that the category of topological $\top$-convergence groups  is isomorphic to the category of topological $\top$-neighborhood groups. In \cite{LL21}, Li and Jin proved that the topological condition in topological $\top$-neighborhood groups is reductant. Inspired by the above two works, we put forward Theorem \ref{themre}.
\end{remark}

\section{The uniformization}

In this section, we will show that  $\top$-convergence groups can be uniformizable. That is, for each $\top$-convergence group $(X,\cdot, \TC)$, there exists a $\top$-uniformly space such that it can reproduces $\TC$.

At first, we fix some notations. For $\lambda,\mu\in L^{X\times X}$, we define $\mu\circ\lambda, \lambda^{\thicksim 1} \in L^{X\times X}$ as $\forall x,y\in X$, $$(\mu\circ\lambda)(x,y)=\bigvee_{z\in X}\Big(\lambda(x,z)\ast \mu(z,y)\Big),\ \ \lambda^{\thicksim 1}(x,y)=\lambda(y,x).$$

 For $\TF,\TG\in \TF_L^\top(X\times X)$, note that the family $\TB=\{\mu\circ \lambda |\lambda\in \TF, \mu\in \TG\}$ forms a $\top$-filter base on $X\times X$ iff $\bigvee\limits_{(x,y)\in X\times X}(\mu\circ \lambda)(x,y)=\top$ for any $\lambda\in \TF, \mu\in \TG$. We denote the $\top$-filter generated by $\TB$ as $\TF\circ \TG$ and say that $\TF\circ \TG$ exists.

For $\TF\in \TF_L^\top(X\times X)$, it is easily proved that the family $\{\lambda^{\thicksim 1}|\lambda\in \TF\}$ forms a $\top$-filter on $X\times X$, and denoted  as $\TF^{\thicksim 1}$.


\begin{definition} A subset $\Phi\subseteq \TF_L^\top(X\times X)$  is called a $\top$-uniformly convergence structure on $X$ if

(TUC1) $\forall x\in X$, $[(x,x)]_{\top}\in \Phi$,

(TUC2) $\TF\subseteq \TG$ and $\TF\in \Phi\Longrightarrow \TG\in \Phi$,


(TUC3) $\TF, \TG \in \Phi$ and $\TF\circ \TG$ exists $\Longrightarrow \TF\circ \TG\in \Phi$,

(TUC4) $\TF\in \Phi\Longrightarrow \TF^{\thicksim 1}\in \Phi$.

The pair $(X,\Phi)$ is called a $\top$-uniformly convergence space.

$(X,\Phi)$ is called a $\top$-uniformly limit space if it fulfills moreover

(TUC4) $\TF,  \TG \in \Phi\Longrightarrow \TF\cap \TG\in \Phi$.

\end{definition}


For a group $(X,\cdot)$ and $\lambda\in L^X$, we define $\lambda_l\in L^{X\times X}$ as $\forall x,y\in X, \lambda_l(x,y)=\lambda(x^{-1}y)$. It is easily seen that $(\lambda_l)^{\thicksim 1}=(\lambda^{-1})_l$. Furthermore, if $\TF\in \TF_L^\top(X)$ then  the family $\TB_l=\{\lambda_l|\lambda\in \TF\}$ forms a $\top$-filter base on $X\times X$, then the generated  $\top$-filter is denoted as $\TF_l$.

The next we  show that each $\top$-convergence group
induces a $\top$-uniformly convergence space. The following lemma is needed.

%
%
%

\begin{lemma} \label{lemfuhe} Let  $\TF,\TG\in \TF_L^\top(X)$ and  $\TF_l\circ \TG_l$ exists. Then $(\TF\odot\TG)_l\subseteq \TF_l\circ \TG_l$.
\end{lemma}

\begin{proof} Let $\lambda\in \TF\odot\TG$, then $\bigvee\limits_{\mu\in \TF,\nu\in \TG} S_X(\mu\odot \nu, \lambda)=\top$. Take any $\mu\in \TF,\nu\in \TG$, then  for any $(x,y)\in X\times X$,
\begin{eqnarray*}(\nu_l\circ \mu_l)(x,y)&=&\bigvee_{z\in X}\Big(\mu_l(x,z)\ast \nu_l(z,y)\Big)\leq\bigvee_{z\in X}\Big(\mu_l(x,z)\wedge \nu_l(z,y)\Big)\\
&=&\bigvee_{z\in X}\Big(\mu(x^{-1}z)\wedge \nu(z^{-1}y)\Big), \ {\rm by}\ (x^{-1}z)\cdot (z^{-1}y)=x^{-1}y\\
&\leq&(\mu\odot \nu)(x^{-1}y)\\
&=&(\mu\odot \nu)_l(x,y).
\end{eqnarray*}It follows that
\begin{eqnarray*}\top&=&\bigvee\limits_{\mu\in \TF,\nu\in \TG} S_X(\mu\odot \nu, \lambda)\\
&=& \bigvee_{\mu\in \TF,\nu\in \TG} \bigwedge_{z\in X}\big[(\mu\odot \nu)(z)\rightarrow \lambda(z)\big]\\
&\leq& \bigvee_{\mu\in \TF,\nu\in \TG} \bigwedge_{(x,y)\in X\times X}\big[(\mu\odot \nu)(x^{-1}y)\rightarrow \lambda(x^{-1}y)\big]\\
&=&\bigvee_{\mu\in \TF,\nu\in \TG} S_{X\times X}((\mu\odot \nu)_l, \lambda_l)\\
&\leq& \bigvee_{\mu\in \TF,\nu\in \TG} S_{X\times X}(\mu_l\circ \nu_l, \lambda_l),\end{eqnarray*} which means $\lambda_l\in \TF_l\circ \TG_l$. Hence $(\TF\odot\TG)_l\subseteq \TF_l\circ \TG_l$.
\end{proof}

\begin{proposition} \label{thmLUl} Let $(X,\cdot, \TC)$ be a $\top$-convergence group. Then the pair $(X,\Phi^{\TC})$ defined by
$$\forall \TF\in \TF_L^{\top}(X\times X), \TF\in \Phi^{\TC}\Longleftrightarrow \exists \TG\stackrel{\TC}{\dashrightarrow} e\ {\rm s.t.}\ \TF\supseteq \TG_l$$is a $\top$-uniformly convergence space.
\end{proposition}

\begin{proof}

(TUC1). Let $x\in X$. For any $\lambda \in [e]_\top$, it follows by $\lambda_l(x,x)=\lambda(x^{-1}x)=\lambda(e)=\top$ that $\lambda_l\in [(x,x)]_\top$, which means $([e]_\top)_l\subseteq [(x,x)]_\top$. From $[e]_\top\stackrel{\TC}{\dashrightarrow} e$ we obtain that $[(x,x)]_\top\in \Phi^{\TC}$.

(TUC2). It is clearly.


(TUC3) Let  $\TF_1, \TF_2 \in \Phi^{\TC}$ and $\TF_1\circ \TF_2$ exists. Then there exist  $\TG_i\in \TF_L^\top(X)$ ($i=1,2$) s.t. $\TG_i\stackrel{\TC}{\dashrightarrow} e$ and $\TF_i\supseteq (\TG_i)_l$. So, $\TG_1\odot \TG_2\stackrel{\TC}{\dashrightarrow} e$ from (TCG1). Note that $\TF_1\circ \TF_2$ exists and $(\TG_1)_l\circ (\TG_2)_l\subseteq \TF_1\circ \TF_2 $, hence $(\TG_1)_l\circ (\TG_2)_l$ exists. It follows by Lemma \ref{lemfuhe} that  $(\TG_1\odot\TG_2)_l\subseteq (\TG_1)_l\circ (\TG_2)_l\subseteq \TF_1\circ \TF_2 $, which means $\TF_1\circ \TF_2\in \Phi^{\TC}$.

(TUC4) Let  $\TF \in \Phi^{\TC}$. Then there exists  $\TG\in \TF_L^\top(X)$ s.t. $\TG\stackrel{\TC}{\dashrightarrow} e$ and $\TF\supseteq \TG_l$. So, $\TG ^{-1}\stackrel{\TC}{\dashrightarrow} e$ from (TCG2). It follows by $\forall \lambda \in L^X$, $(\lambda_l)^{\thicksim 1}=(\lambda^{-1})_l$ that  $(\TG^{-1})_l=(\TG_l)^{\thicksim 1}\subseteq \TF^{\thicksim 1}$. Hence $\TF^{\thicksim 1}\in \Phi^{\TC}$.
\end{proof}

Next we exhibit  that each $\top$-uniformly convergence space
induces a $\top$-convergence space. The following lemma is needed.

\begin{lemma} \label{lemmatoc} Let $\TF, \TG\in \TF_L^\top(X)$ and $x\in X$. Then $\TG_l\subseteq [x]_{\top}\times \TF\Longleftrightarrow [x]_\top\odot \TG\subseteq \TF$.
\end{lemma}

\begin{proof} $\Longrightarrow$. Assume that $\TG_l\subseteq [x]_{\top}\times \TF$. Take $\lambda\in [x]_\top, \mu\in
\TG$, then  $\lambda\odot \mu\geq
\top_{\{x\}}\odot \mu=\mu_l(x,-)$. From $\mu\in
\TG$ and $\TG_l\subseteq [x]_{\top}\times \TF$, we know that $\mu_l\in [x]_{\top}\times \TF$, so $$\top=\bigvee_{\nu\in \TF}S_{X\times X}(\top_{\{x\}}\times \nu, \mu_l)=\bigvee_{\nu\in \TF}S_{X}(\nu, \mu_l(x,-)),$$which means $\mu_l(x,-)\in \TF$. It follows by $\lambda\odot \mu\geq\mu_l(x,-)$ we get $\lambda\odot \mu\in \TF$. Therefore, $[x]_\top\odot \TG\subseteq \TF$.

$\Longleftarrow$. Assume that $[x]_\top\odot \TG\subseteq \TF$. Take $\lambda\in
\TG$, then $\top_{\{x\}}\odot \lambda \in [x]_\top\odot \TG$. Notice that $\top_{\{x\}}\odot \lambda=\lambda_l(x,-)$, hence
$\lambda_l(x,-)\in [x]_\top\odot \TG\subseteq\TF$, and so $\top_{\{x\}}\times \lambda_l(x,-)\in [x]_{\top}\times \TF$.  Then from $\top_{\{x\}}\times \lambda_l(x,-)\leq \lambda_l$ one get $\lambda_l\in  [x]_{\top}\times \TF$. Therefore, $\TG_l\subseteq [x]_{\top}\times \TF$.
\end{proof}

\begin{proposition} \label{thmUL1} Let $(X,\Phi)$ be a $\top$-uniformly convergence space. Then the pair $(X,\TC_\Phi)$ defined by $$\forall \TF\in \TF_L^\top(X), \forall x\in X,  \TF\stackrel{\TC_\Phi}{\dashrightarrow} x\Longleftrightarrow [x]_{\top}\times \TF\in \Phi$$
is a $\top$-convergence space.
\end{proposition}

\begin{proof} (TC1) For any $x\in X$, it follows by  $[x]_\top\times [x]_\top=[(x,x)]_\top\in \Phi$, we get $[x]_\top\stackrel{\TC_\Phi}{\dashrightarrow} x$.

(TC2) It is clearly.
\end{proof}

\begin{definition} A  $\top$-convergence space $(X, \TC)$  is named $\top$-uniformizable if $\TC=\TC_{\Phi}$ for some $\top$-uniformly convergence space $(X,\Phi)$.
\end{definition}

The next theorem shows that each  $\top$-convergence group is $\top$-uniformizable.

\begin{theorem} (Uniformizable theorem) For a $\top$-convergence group $(X,\cdot, \TC)$, we have $\TC=\TC_{\Phi_{\TC}}$.
\end{theorem}

\begin{proof} Let $ \TF\stackrel{\TC_{\Phi_{\TC}}}{\dashrightarrow} x$. Then $[x]_\top\times \TF\in \Phi_{\TC}$, so there is a $\TG\in \TF_L^\top(X)$ s.t. $\TG\stackrel{\TC}{\dashrightarrow} e$ and $\TG_l\subseteq [x]_\top\times \TF$. Notice that
$$\TG\stackrel{\TC}{\dashrightarrow} e, [x]_\top \stackrel{\TC}{\dashrightarrow} x \Longrightarrow [x]_\top\odot \TG\stackrel{\TC}{\dashrightarrow} x\ {\rm and}\ \TG_l\subseteq [x]_\top\times \TF\stackrel{\rm Lemma \ref{lemmatoc}}{\Longrightarrow} [x]_\top\odot \TG\subseteq \TF,$$which means $\TF\stackrel{\TC}{\dashrightarrow} x$.

Conversely, let $\TF\stackrel{\TC}{\dashrightarrow} x$, then $\TG:=[x^{-1}]_\top\odot \TF\stackrel{\TC}{\dashrightarrow} e$ by Theorem \ref{theorem01}. It follows by $[x]_\top\odot \TG=\TF$ and Lemma \ref{lemmatoc}  one get $\TG_l\subseteq [x]_\top\times \TF$, and so $[x]_\top\times \TF\in \Phi_{\TC}$ by $\TG\stackrel{\TC}{\dashrightarrow} e$. Hence $ \TF\stackrel{\TC_{\Phi_{\TC}}}{\dashrightarrow} x$ as desired.
\end{proof}


\section{The power object}

In \cite{Yu17}, Yu and Fang proved that the category {\bf TCons}  has power object when $L$ is restricted to be  a complete MV-algebra. In this section, we will further show that the category {\bf TConG} also has power object when the underlying lattice of $L$ satisfying a distributive law (CD)\footnote{Many complete residuated lattices satisfy  (CD) law, for example, BL-algebra, MV-algebra, $\prod$-algebra, Heyting algebra, left continuous t-norm, etc.}: $$\forall a\in L, \forall \{b_j\}(j\in J)\in L, a\wedge\bigvee\limits_{j\in J} b_j=\bigvee\limits_{j\in J}(a\wedge b_j).$$

Initially, we summarize Yu and Fang's findings in \cite{Yu17}. We observe that the majority of their results are applicable to an arbitrary complete residuated lattice $L$. However, in certain instances, the underlying lattice $L$ must adhere to the (CD) law. Consequently, throughout this section, we presuppose that $L$ satisfies the (CD) law.

\begin{lemma} \label{lemcd} (\cite{Yu17}) Let $f_i:X_i\longrightarrow Y_i$ be mappings and $\TF_i\in \TF_L^\top(X_i)$ ($i=\{1,2\}$). Then $(f_1\times f_2)^\Rightarrow(\TF_1\times \TF_2)= f_1^\Rightarrow(\TF_1) \times  f_2^\Rightarrow(\TF_2)$.
\end{lemma}

\begin{proof} By using the (CD) law we can verify that $(f_1\times f_2)^\rightarrow (\lambda_1\times \lambda_2)=f_1^\rightarrow (\lambda_1) \times f_2^\rightarrow(\lambda_2)$ for any $\lambda_i\in \TF_L^\top(X_i)$. Therefore, $(f_1\times f_2)^\Rightarrow(\TF_1\times \TF_2)$ and $ f_1^\Rightarrow(\TF_1) \times  f_2^\Rightarrow(\TF_2)$ have a common $\top$-filter base $\{f_1^\rightarrow (\lambda_1) \times f_2^\rightarrow(\lambda_2)|\lambda_i\in \TF_L^\top(X_i)\}$, and hence they are equal.
\end{proof}

\begin{definition} (\cite{Yu17}) \label{defnsh} Let $(X,\TC_X), (Y,\TC_Y)\in $ {\bf TCons} and $$C(X,Y)=\{f:(X,\TC_X)\longrightarrow (Y,\TC_Y)|f \ {\rm is \  continuous}\}.$$ Then the mapping $ev: C(X,Y)\times X\longrightarrow Y$ defined by $$\forall f\in C(X,Y), \forall x\in X, ev (f,x)=f(x),$$ is called the evaluation mapping.
Moreover, the mapping $\TC: \TF_L^\top(C(X,Y))\longrightarrow 2^{C(X,Y)}$ defined by $\forall \TF\in \TF_L^\top(C(X,Y)),\forall f\in C(X,Y)$,
$$\TF\stackrel{\TC}{\dashrightarrow} f\Longleftrightarrow \forall x\in X, \forall \TG\in \TF_L^\top(X), \TG\stackrel{\TC_{X}}{\dashrightarrow} x \ {\rm implies}\ ev^\Rightarrow(\TF\times \TG)\stackrel{{\TC_{Y}}}{\dashrightarrow} f(x),$$ is a $\top$-convergence structure on $C(X,Y)$, called the power object in {\bf TCons}.\end{definition}

\begin{proposition} \label{funs} Let $(X,\TC_X)$, $(Y,\TC_Y)$, $(Z,\TC_Z) \in $ {\bf TCons}.

{\rm (1)} The evaluation mapping $ev$ is continuous (\cite[Proposition 3.8]{Yu17}).

{\rm (2)} If the mapping $\varphi: (Z\times X, \TC_Z\times \TC_X)\longrightarrow (Y,\TC_Y)$ is continuous then  the mapping $$\varphi^\lozenge: (Z,\TC_Z)\longrightarrow (C(X,Y), \TC): z\longmapsto \varphi(z,-)$$ is the unique continuous mapping s.t. $ev\circ (\varphi^\lozenge\times id_X)=\varphi$.  (\cite[Lemma 3.10]{Yu17}). In that, Lemma \ref{lemcd} is used.
\end{proposition}

\begin{remark} It follows by Proposition \ref{funs}  that the category {\bf TCons} is Cartesian closed, see \cite[Theorem 3.11]{Yu17}.\end{remark}

%
%
Note that by replacing the continuous mappings in {\bf TCons} with continuous group homomorphisms, one can define the notions of source and initial structure in {\bf TConG}.

Let $\big((X,\cdot) \stackrel{f_j}{\longrightarrow} (X_j,\cdot,\TC_j)\big)_{j\in J}$ be a source in {\bf TConG}. Zhang and Pang \cite{LZ23} showed that the triple $(X,\cdot, \TC)$ forms a $\top$-convergence group, where $(X,\TC)$ is the initial structure w.r.t the source $\big(X \stackrel{f_j}{\longrightarrow} (X_j, \TC_j)\big)_{j\in J}$ in {\bf TCons}. Therefore, they obtained the following theorem.

\begin{theorem} (\cite{LZ23}) \label{tngin} Each source $\big((X,\cdot) \stackrel{f_j}{\longrightarrow} (X_j,\cdot,\TC_j)\big)_{j\in J}$ in {\bf TConG} has an initial structure. That is, {\bf TConG} is  topological over {\bf Grp} w.r.t.
the forgetful functor.
\end{theorem}

Because every topological category has  products, the category {\bf TConG} has products, and therefore has finite products. Moreover, after forgetting the group operations, the product of the objects in {\bf TConG} happens to be their product in {\bf TCons}.

Next, we turn our attention to the power object of the category {\bf TConG}. Some preparation is needed.

\begin{lemma} \label{unv} Let $(X,\TC_X)$, $(Y,\TC_Y)$, $(Z,\TC_Z) \in $ {\bf TCons} and $\varphi: (Z,\TC_Z)\longrightarrow (C(X,Y), \TC)$ be a mapping. If the mapping $$\varphi^\Box: (Z\times X, \TC_Z\times \TC_X)\longrightarrow (Y,\TC_Y), (z,x)\longmapsto \varphi(z)(x)$$ is continuous, then so is $\varphi$.
\end{lemma}

\begin{proof} It follows by Proposition \ref{funs} (2) and $\varphi=\varphi^{\Box\lozenge}$. Indeed, $\forall z\in Z$, $\varphi^{\Box\lozenge}(z)=\varphi^{\Box}(z,-)=\varphi(z).$
\end{proof}

\begin{proposition} \label{Cgu1}  Let $(X,\cdot_X, \TC_X), (Y,\cdot_Y, \TC_Y)\in $ {\bf TConG} and $C(X,Y)$ be the set defined  in Definition  \ref{defnsh} via $(X, \TC_X), (Y, \TC_Y)\in $ {\bf TCons}. Then the pair $(C(X,Y), \cdot)$ forms  a group, where the group operation $\cdot$ is defined as $$\forall f,g\in C(X,Y), \forall x\in X,  (f\cdot g)(x)=f(x)\cdot_Y g(x).$$For the convenience of expression, we also write the above group operation $\cdot$ as $m$.
\end{proposition}

\begin{proof} (1) For any $f,g\in C(X,Y)$, we check below $f\cdot g\in C(X,Y)$.

Let $\TF\in \TF_L^\top(X)$ and $\TF\stackrel{\TC_X}{\dashrightarrow} x$. We prove below $f^\Rightarrow(\TF)\odot g^\Rightarrow(\TF)\subseteq (f\cdot g)^\Rightarrow(\TF)$. Indeed, notice that $$\TB=\{\lambda \odot \mu|\lambda\in f^\Rightarrow(\TF), \mu\in g^\Rightarrow(\TF)\}$$ forms a $\top$-filter base of  $f^\Rightarrow(\TF)\odot g^\Rightarrow(\TF)$. Put $\lambda\in f^\Rightarrow(\TF), \mu\in g^\Rightarrow(\TF)$, then $f^\leftarrow(\lambda), g^\leftarrow(\mu) \in \TF$, and so $f^\leftarrow(\lambda)\wedge g^\leftarrow(\mu) \in \TF$. For any $y\in Y$,
\begin{eqnarray*}(f\cdot g)^\rightarrow \Big(f^\leftarrow(\lambda)\wedge g^\leftarrow(\mu)\Big)(y)&=&\bigvee_{(f\cdot g)(x)=y}\Big(f^\leftarrow(\lambda)(x)\wedge g^\leftarrow(\mu)(x)\Big)\\
&=&\bigvee_{f(x)\cdot_Y g(x)=y}\Big(\lambda (f(x))\wedge \mu(g(x))\Big)\\
&\leq&\bigvee_{y_1\cdot_Y y_2=y}\Big(\lambda (y_1)\wedge \mu(y_2)\Big)\\
&=&(\lambda \odot \mu )(y),
\end{eqnarray*} i.e., $(f\cdot g)^\rightarrow \Big(f^\leftarrow(\lambda)\wedge g^\leftarrow(\mu)\Big)\leq \lambda \odot \mu$. It follows by $f^\leftarrow(\lambda)\wedge g^\leftarrow(\lambda) \in \TF$ we get $\lambda \odot \mu\in (f\cdot g)^\Rightarrow(\TF)$. Hence $f^\Rightarrow(\TF)\odot g^\Rightarrow(\TF)\subseteq (f\cdot g)^\Rightarrow(\TF)$.

By $f,g\in C(X,Y)$ and $\TF\stackrel{\TC_X}{\dashrightarrow} x$, we have $f^\Rightarrow(\TF)\stackrel{\TC_Y}{\dashrightarrow} f(x)$ and $g^\Rightarrow(\TF)\stackrel{\TC_Y}{\dashrightarrow} g(x)$, it follows by $(Y,\cdot_Y, \TC_Y)\in $ {\bf TConG} that $f^\Rightarrow(\TF)\odot g^\Rightarrow(\TF)\stackrel{\TC_Y}{\dashrightarrow} (f\cdot g)(x)$, and so $(f\cdot g)^\Rightarrow(\TF)\stackrel{\TC_Y}{\dashrightarrow} (f\cdot g)(x)$ since $f^\Rightarrow(\TF)\odot g^\Rightarrow(\TF)\subseteq (f\cdot g)^\Rightarrow(\TF)$. Now, we have proved that $f\cdot g\in C(X,Y)$.

(2) Obviously, the operation $\cdot$ on $C(X,Y)$ satisfies the associative law.

(3) Let $e_Y$ be the identity element of the group $(Y,\cdot_Y)$. We  use $\overline{e_Y}$ to denote the constant-value mapping from $X$ to $Y$ values $e_Y$. Then it is easily seen that $\overline{e_Y}\in C(X,Y)$ and $f\cdot \overline{e_Y}=\overline{e_Y}\cdot f=f$ for any $f\in C(X,Y)$. Hence $\overline{e_Y}$ is  the identity element of $(C(X,Y),\cdot)$.

(4) For any $f\in C(X,Y)$, we define $r(f)=r_Y\circ f$. Then $r(f)\in C(X,Y)$ since $r_Y$ and $f$ are all continuous mappings. In addition,  $f\cdot r(f)=\overline{e_Y}=r(f)\cdot f$.

From (1)-(4) we know that $(C(X,Y), \cdot)$ is a group with the identity $\overline{e_Y}$ and the inverse operation $r$.
\end{proof}

\begin{theorem} \label{Cgu}  Let $(X,\cdot_X, \TC_X), (Y,\cdot_Y, \TC_Y)\in $ {\bf TConG} and $(C(X,Y), \TC)$ be the $\top$-convergence space defined  in Definition  \ref{defnsh}. Then the triple $(C(X,Y), \cdot, \TC)\in $ {\bf TConG}.

\end{theorem}

\begin{proof}  We only need to prove that the group operations $m$ and $r$ are continuous.

(1) We prove below that $m:C(X,Y)\times C(X,Y)\longrightarrow C(X,Y)$ is continuous.

From Lemma  \ref{unv}, we need  only check that $m^\Box:C(X,Y)\times C(X,Y)\times X\longrightarrow Y$ is continuous. We define a mapping $$\psi:C(X,Y)\times C(X,Y)\times X\longrightarrow  \Big(C(X,Y)\times X\Big) \times \Big( C(X,Y)\times X\Big), (f,g,x)\longmapsto \Big((f,x), (g,x)\Big).$$ It follows by the initial property of product spaces that the mapping $\psi$ is continuous. In addition, from Proposition \ref{conpr} we know that the mapping $ev\times ev: \Big(C(X,Y)\times X\Big) \times \Big( C(X,Y)\times X\Big)\longrightarrow Y\times Y$ is continuous. Obviously, the mapping $m_Y: Y\times Y\longrightarrow Y$ is continuous. We prove below that  $m^\Box=m_Y\circ (ev\times ev)\circ \psi$, i.e., the following diagram commutes. $$\bfig \morphism(0,0)|a|/@{->}@<0pt>/<1800,0>[C(X,Y)\times C(X,Y)\times X`
\Big(C(X,Y)\times X\Big) \times \Big( C(X,Y)\times X\Big);\psi] \morphism(0,0)|l|/@{->}@<0pt>/<0,-400>[C(X,Y)\times C(X,Y)\times X`Y;m^\Box]
\morphism(1800,0)|r|/@{->}@<0pt>/<0,-400>[\Big(C(X,Y)\times X\Big) \times \Big( C(X,Y)\times X\Big)`Y\times Y;ev\times ev] \morphism(0,-400)|a|/@{<-}@<0pt>/<1800,0>[Y`Y\times Y;m_Y] \efig
$$
Indeed, $\forall (f,g,x)\in C(X,Y)\times C(X,Y)\times X$,
$$m_Y\circ (ev\times ev)\circ \psi(f,g,x)=ev(f,x)\cdot_Y ev(g,x)=f(x)\cdot_Y g(x)=(f\cdot g)(x)=m(f,g)(x)=m^\Box(f,g,x).$$

Hence $m^\Box=m_Y\circ (ev\times ev)\circ \psi$ is continuous, as desired.

(2) We prove below that $r:C(X,Y)\longrightarrow C(X,Y), f\longmapsto r_Y\circ f$ is continuous.

From Lemma  \ref{unv}, we need  only check that $r^\Box:C(X,Y)\times X\longrightarrow Y$ is continuous. Notice that $$\forall (f,x)\in C(X,Y)\times X, r^\Box(f,x)=r(f)(x)=r_Y( f(x))=(r_Y\circ ev) (f,x),$$ which means that $r^\Box=r_Y\circ ev$ is continuous, as desired.

Now, we have completed the proof of that $(C(X,Y), \cdot, \TC)$ is a $\top$-convergence group.
\end{proof}

\begin{definition} Let $(X,\cdot_X, \TC_X), (Y,\cdot_Y, \TC_Y)\in $ {\bf TConG}. Then the $\top$-convergence group $(C(X,Y), \cdot, \TC)$ defined in Theorem  \ref{Cgu} is called the power object w.r.t. $(X,\cdot_X, \TC_X)$ and $(Y,\cdot_Y, \TC_Y)$.

\end{definition}

\begin{remark} Note that the evaluation mapping $ev$ is not a group
homomorphism generally. Hence, we can not prove as Yu and Fang \cite{Yu17} that {\bf TConG} is Cartesian closed.

\end{remark}

%
%
%
%
%

\section{Concluding remarks}
In this paper, we have conducted an in-depth investigation into the category of $\top$-convergence groups. We introduced a fundamental characterization of $\top$-convergence groups through the $\odot$-product of $\top$-filters, and substantiated that this category exhibits localization, uniformization, and power object properties. Moreover, we demonstrated that the topological condition (TT) in the category of topological $\top$-convergence groups is superfluous. As highlighted in the introduction, there exist similar problems within enriched $L$-generalized convergence groups that remain unexplored. Consequently, these unexplored problems will be addressed in our subsequent research endeavors. Additionally, various types of lattice-valued convergence spaces, defined by fuzzy nets \cite{XY12} and fuzzy ideals \cite{MG05}, exist. We intend to study lattice-valued convergence groups based on these spaces in the future.


\begin{thebibliography}{10}



\bibitem{AHS} J. Ad\'{a}mek, H. Herrlich, G.E. Strecker,  Abstract and
Concrete Categories, Wiley, New York, 1990.



\bibitem{TA88}  T.M.G. Ahsanullah, On fuzzy neighborhood groups, Journal of Mathematical Analysis
and Applications, 130 (1988) 237--251.

\bibitem{TA08}  T.M.G. Ahsanullah, J. Al-Mufarrij, Framed-valued stratified
generalized convergence groups, Quaestiones Mathematicae, 31
(2008) 279--302.



\bibitem{TG142}  T.M.G. Ahsanullah, D. Gauld, J. Al-Mufarrij, F. Al-Thukair, Enriched lattice-valued convergence groups, Fuzzy Sets and Systems, 238 (2014) 71--88.


\bibitem{TA17}  T.M.G. Ahsanullah, G. J\"{a}ger, Stratified $LMN$-convergence tower groups and their stratified
$LMN$-uniform convergence tower structures, Fuzzy Sets and Systems, 330 (2018) 105--123.

\bibitem{TA19} T.M.G. Ahsanullah, G. J\"{a}ger, Quantale-valued generalizations of approach groups, New Mathematics and Natural Computation, 15,1 (2019) 1--30.

\bibitem{JA08} J. Al-Mufarrij, T. M. G. Ahsanullah, On the category of fixed
basis frame valued topological groups, Fuzzy Sets and Systems, 159 (2008) 2529--2551.

\bibitem{AA08} A. Arhangel'skii, M. Tkachenko, Topological groups
and related structures, Atlantis Press World Scientific, Paris, 2008.

\bibitem{FB05} F. Bayoumi, On initial and final $L$-topological groups, Fuzzy Sets
and Systems, 156 (2005) 43--54.

\bibitem {FB08} F. Bayoumi, Global $L$-neighborhood groups, Fuzzy Sets and
Systems, 159 (2008) 605--619.



\bibitem{R.B02} R. B\v{e}lohl\'{a}vek,  Fuzzy Relational Systems:
Foundations and Principles, Kluwer Academic Publishers, New York,
2002.





\bibitem{JF04} J.M. Fang, $I$-FTOP is isomorphic to $I$-FQN and I-AITOP, Fuzzy Sets
and Systems, 147 (2004) 317--325.


\bibitem{J.F10} J.M. Fang, Stratified $L$-ordered convergence structures,
Fuzzy Sets and Systems, 161 (2010) 2130--2149.

\bibitem{fang2017} J.M. Fang, Y.L. Yue, $\top$-diagonal conditions and Continuous extension theorem, Fuzzy Sets and Systems, 321 (2017) 73--89.



\bibitem{JF21} J.M. Fang, Y.L. Yue, Extensionality and E-connectedness in the category of $\top$-convergence spaces, Fuzzy Sets and Systems, 425 (2021) 100--116.

\bibitem{JF23} J.M. Fang, Y.L. Yue, A one to one correspondence between connected varieties and disconnected varieties of $\top$-semiuniform convergence spaces,
Fuzzy Sets and Systems,  466 (2023) 108418.

\bibitem{DF79} D.H. Foster, Fuzzy topological
groups, Journal of Mathematical Analysis and Applications,
 67 (1979) 549--564.




\bibitem{MG05} M. G\"{u}lo\u{g}lu, D. \c{C}oler, Convergence in I-fuzzy topological spaces, Fuzzy Sets and Systems, 151 (2005) 615--623.



\bibitem{H98} P. H\'{a}jek, Metamathematics of Fuzzy Logic,  Kluwer
Academic Publishers, Dordrecht, 1998.












\bibitem{U.H99} U. H\"{o}hle, A. \u{S}ostak, Axiomatic foundations of
fixed-basis fuzzy topology, in: U. H\"{o}hle, S.E. Rodabaugh (Eds.),
Mathematics of Fuzzy Sets: Logic, Topology and Measure Theory, The
Handbooks of Fuzzy Sets Series, Vol.3, Kluwer Academic Publishers,
Boston, Dordrecht, London, 1999, pp.123--273.

\bibitem{QL18} X.K. Huang, Q.G. Li, Q.M. Xiao, The $L$-ordered semigroups based on $L$-partial orders, Fuzzy Sets and Systems, 339 (2018) 31--50.










\bibitem{L.L171} Q. Jin,  L.Q. Li, Stratified lattice-valued neighborhood tower group, Quaestiones Mathematicae, 41,6 (2018) 847-861.







\bibitem{G.J01} G. J\"{a}ger, A category of $L$-fuzzy convergence spaces,
Quaestiones Mathematicae, 24 (2001) 501--517.



\bibitem{G.J07} G. J\"{a}ger, Pretopological and topological lattice-valued convergence spaces,  Fuzzy Sets and Systems, 158 (2007) 424--435.




\bibitem{G.J15} G. J\"{a}ger, Stratified $LMN$-convergence tower spaces,  Fuzzy Sets and Systems,  282 (2016) 62--73.

\bibitem{GJ19} G. J\"{a}ger,
Quantale-valued generalizations of approach spaces and quantale-valued topological spaces,  Quaestiones Mathematicae, 42,10 (2019) 1313--1333.





\bibitem{G.J24} G. J\"{a}ger, Uniformly continuous extension in $L$-uniform convergence tower spaces,  Filomat,  38,2 (2024) 577--588.

\bibitem{YY22} G. J\"{a}ger, Y.L. Yue $\top$-uniform convergence spaces, Iranian Journal of Fuzzy Systems, 19, 2 (2022) 133--149.

\bibitem{HL17} H.L. Lai, D.X. Zhang,  Fuzzy topological spaces with conical neighborhood
system,  Fuzzy Sets and Systems, 330 (2018) 87--104.


%





\bibitem{LL12} L.Q. Li, Q. Jin, On stratified $L$-convergence
spaces: Pretopological axioms and diagonal axioms, Fuzzy Sets and
Systems, 204 (2012) 40--52.

\bibitem{LL21} L.Q. Li, Q. Jin, A category of complete residuated lattice-value neighborhood groups, Fuzzy Sets and
Systems, 442 (2022) 53--75.

\bibitem{LL20} L.Q. Li, Q. Jin,  C.X. Bo, Z.Y. Xiu, The categorical relationships between neighborhood spaces, $\top$-neighborhood spaces and stratified $L$-neighborhood spaces, Filomat, 35,4 (2021) 1267--1287.



\bibitem{LL18} L.Q. Li, Q. Jin, K. Hu, Lattice-valued convergence associated with
CNS spaces, Fuzzy Sets and Systems, 370(2019) 91-98.







\bibitem{RL76} R. Lowen,  Fuzzy topological spaces and fuzzy compactness, Journal of Mathematical analysis and applications, 56 (1976) 621--633.




\bibitem{RL82} R. Lowen, Fuzzy neighborhood spaces, Fuzzy Sets and Systems,  7(1982) 165--189.





\bibitem{P17} B. Pang, Y. Zhao,  Several types of enriched $(L, M)$-fuzzy convergence spaces, Fuzzy Sets and Systems, 321 (2017) 55--72.

\bibitem{P172} B. Pang,  Strafified $L$-ordered filter spaces, Quaestiones Mathematicae, 40,5 (2017) 661--678.

\bibitem{P18} B. Pang, Z.Y. Xiu, Stratified $L$-prefilter convergence structures in stratified $L$-topological spaces,  Soft Computing, 22,22 (2018) 7539--7551.

%




\bibitem{GR} L. Reid, G. Richardson, Connecting $\top$ and Lattice-Valued
Convergences, Iranian Journal of Fuzzy Systems, 15,4 (2018) 151--169.

\bibitem{GR18} L. Reid, G. Richardson, Lattice-valued spaces: $\top$-Completions, Fuzzy Sets and
Syststems, 369(2019) 1--19.



%

\bibitem{SH92} J.Z. Shen, Fuzzifying topological groups based on completely distributive residuated lattice-valued logic (I), International symposium on multiple-valued logic, (1992) 198--205.

\bibitem{SH94} J.Z. Shen, Fuzzifying topological groups based on completely distributive residuated lattice-valued logic (II), Information Sciences, 80 (1994) 319--339.

\bibitem{FS09} F.G. Shi, $L$-fuzzy interiors and $L$-fuzzy closures, Fuzzy Sets and Syststems, 160 (2009) 1218--1232.









\bibitem{CY10} C.H. Yan, S.Z. Guo, $I$-fuzzy topological groups, Fuzzy Sets
and Systems, 161 (2010) 2166--2180.

\bibitem{XY12} X.F. Yang, S.G. Li, Net-theoretical convergence in $(L,M)$-fuzzy cotopologicalspaces, Fuzzy Sets and Systems, 204 (2012) 53--65.



\bibitem{WY12} W. Yao, On many-valued stratified $L$-fuzzy convergence spaces, Fuzzy Sets and Systems, 159 (2008) 2503--2519.


\bibitem{Yu17} Q. Yu, J.M. Fang, The Category of $\top$-convergence spaces
and its cartesian-closedness, Iranian Journal of Fuzzy Systems, 14,3
(2017) 121--138.

\bibitem{Yue20} Y.L. Yue, J.M. Fang, The $\top$-monad and its applications, Fuzzy Sets and Systems, 382 (2020) 79--97.

\bibitem{Yue21} Y.L. Yue, J.M. Fang, W. Yao,  Monadic convergence structures revisited, Fuzzy Sets and Systems, 406 (2021) 107--118.


\bibitem{DZ07} D.X. Zhang, An enriched category approach to many valued topology, Fuzzy Sets and Systems, 158 (2007) 349--366.

\bibitem{LZ22} L. Zhang, B. Pang, Monoidal closedness of the category of $\top$-semiuniform convergence spaces,  Hacettepe Journal of Mathematics and Statistics, 51, 5 (2022) 1348--1370.

\bibitem{LZ23} L. Zhang, B. Pang, A new approach to lattice-valued convergence groups via $\top$-filters,  Hacettepe Journal of Mathematics and Statistics, Fuzzy Sets and Systems, 455 (2023) 198--221.

\bibitem{Yan12} S.Y. Zhang, C.H. Yan, $L$-fuzzifying topological groups, Iranian Journal of Fuzzy Systems, 4 (2012) 115--132.





\bibitem{DZ14} H. Zhao, S.G. Li, G.X. Chen, $(L,M)$-fuzzy topological groups, Journal of Intelligent and Fuzzy Systems, 26 (2014) 1517--1526.



\bibitem{HZ07}  H.P. Zhang, J.X. Fang,  $I(L)$-topological groups and its level $L$-topological groups, Fuzzy Sets and Systems, 158 (2007) 1504--1510.
\end{thebibliography}
\end{document}